\documentclass[eqsecnum,preprint,showpacs,amsmath,aps,nofootinbib]{revtex4} 
\usepackage{amsmath}
\usepackage{amssymb}

\begin{document}
\title{From pseudoholomorphic functions to the associated real manifold}

\author{Giampiero Esposito}
\email[E-mail: ]{gesposit@na.infn.it}
\affiliation{Istituto Nazionale di Fisica Nucleare, Sezione di
Napoli, Complesso Universitario di Monte S. Angelo, 
Via Cintia Edificio 6, 80126 Napoli, Italy}

\author{Raju Roychowdhury} 
\email[E-mail: ]{raju@if.usp.br}
\affiliation{Instituto de Fisica, Universidade de Sao Paulo, C. Postal 66318,
05314-970 Sao Paulo, Brazil}

\date{\today}

\begin{abstract}
This paper studies first the differential inequalities that 
make it possible to build a global theory of pseudoholomorphic functions
in the case of one or several complex variables. In the case of one complex dimension, we prove
that the differential inequalities describing pseudoholomorphicity can be used to define a
one-real-dimensional manifold (by the vanishing of a function with nonzero gradient), 
which is here a $1$-parameter family of plane curves. On studying the associated envelopes,
such a parameter can be eliminated by solving two nonlinear partial differential equations.
The classical differential geometry of curves can be therefore exploited to get a novel perspective
on the equations describing the global theory of pseudoholomorphic functions.
\end{abstract}


\maketitle

\section{Introduction}

The progress in complex analysis and differential geometry has led to many important concepts
in pure mathematics and mathematical physics from the nineteenth century until recent times.
For example, when twistor theory and its applications to 
general relativity were developed by Penrose \cite{Penrose1975}, 
and his school, the subject of complex general relativity 
\cite{gesposit-complex} emerged as a fascinating branch of
modern mathematical physics, where the tools of complex differential geometry were applied
to find self-dual or anti-self-dual solutions of (vacuum) Einstein equations
\cite{Penrose1976}, and also to develop
a suitable definition of twistor in curved spacetime \cite{Penrose1975}, without relying on the differential 
equation that defines Killing spinor fields, but rather considering suitable surfaces (e.g. the
totally null $\alpha$- and $\beta$-surfaces \cite{Penrose1986}).

The very concept of complex manifold involves complex-analytic transition functions. More
precisely, a complex manifold is meant to be a paracompact Hausdorff space covered by
neighbourhoods each homeomorphic to an open set in $C^{m}$ ($m=1,2,...$), such that, where these
neighbourhoods overlap, the local coordinates transform by a holomorphic transformation.
Thus, if $z^{1},...,z^{m}$ are complex local coordinates in one such neighbourhood, and if
$w^{1},...,w^{m}$ are local coordinates in another neighbourhood, where they are both defined
one has $w^{i}=w^{i}(z^{1},...,z^{m})$ and each $w^{i}$ is a holomorphic function (see below)
of the $z's$, and the determinant $\partial(w^{1},...,w^{m})/\partial(z^{1},...,z^{m})$ does not
vanish. Well known examples of this abstract concept are the space $C^{m}$, complex projective 
space $CP^{m}$, non-singular submanifolds of $CP^{m}$, the complex torus, 
the orientable surfaces \cite{Chern1979}.

From the point of view of complex analysis, the assumption of holomorphic transition functions is
nontrivial. By definition, the function
$$
f:z=x+{\rm i}y \rightarrow f(z)=u(x,y)+{\rm i}v(x,y)
$$
is holomorphic if it is a continuous function of the complex variable $z$, for which the first
derivative $f'(z)$ exists. This is enough to ensure continuity of $f'(z)$ 
as well \cite{Goursat1884,Goursat1900},
jointly with the many properties that one learns in introductory courses, including the equivalence
with the Weierstrass definition of complex-analytic function of a complex variable, that involves
absolute and uniform convergence of a power series in the first place. However, if the assumption of
differentiability is no longer made, one can define a one-parameter family of functions of a complex
variable, that are holomorphic only if the parameter $\mu$ in the definition is set to $1$. These are the
pseudo-holomorphic functions, beautifully presented by Bers \cite{Bers1956}, but we here rely upon 
the Caccioppoli approach \cite{Caccioppoli1952,Caccioppoli1953}, better suited 
if one wants to build a global theory, but apparently (much) less known
in the literature, maybe because of lack of an English translation. Interestingly, in the theory of
pseudoanalytic functions, no use is made of analytic functions' theory, but on the contrary important
topics of the latter, e.g. the Picard theorem, appear in a new light, through the analysis of
qualitative aspects. In other words, one can obtain elementary proofs of classical theorems,
revealing their intimate nature of metric and topological properties. The Cauchy-Riemann 
differential equations 
$$
u_{x}=v_{y}, \; u_{y}=-v_{x}
$$
that express the holomorphic nature of a function are then replaced, in the pseudoholomorphic case, by
differential inequalities, i.e. simple majorizations replace the equalities among angles of the
conformal representation. 

Section II defines pseudoholomorphic functions of a complex variable according to Caccioppoli, studying in
detail some basic equations in such a definition. Section III presents our definition of pseudoholomorphic 
function of several (i.e. two or more) complex variables, inspired by the Caccioppoli work in the case
of a single complex variable. Since such functions are believed to be more fundamental, one may hope that
the manifolds one arrives at by means of them are also more fundamental in a suitable sense. 
After a review of discwise quasi-conformal maps with $n$ complex variables in Sect. $4$, we study
in Sect. $5$ the $1$-parameter family of plane curves associated with a pseudoholomorphic 
function, while the envelopes for such curves are considered in Sect. $6$. 
Concluding remarks and open problems are presented in Sect. $7$, while some
technical points are discussed in the appendix.

\section{Pseudoholomorphic functions of a complex variable}

Let $z$ be the familiar notation for complex variable $z=x+{\rm i}y$, $x \in {\bf R}, y \in
{\bf R}$, and let $w$ be the continuous function with image
\begin{equation}
w(z)=u(x,y)+{\rm i}v(x,y)
\label{(2.1)}
\end{equation}
defined in a field $A$, a bounded open set of the $z$ plane such that all its points are internal points.
Let the functions $u(x,y),v(x,y)$ satisfy the following assumptions:
\vskip 0.3cm
\noindent
(i) $u(x,y)$ and $v(x,y)$ are absolutely continuous in $x$ and $y$ for almost all values of $y$
and $x$ respectively, while their first derivatives $u_{x},u_{y},v_{x},v_{y}$ are square-integrable
in every internal portion of $A$.
\vskip 0.3cm
\noindent
(ii) If
\begin{equation}
J={\partial (u,v) \over \partial (x,y)}=u_{x}v_{y}-u_{y}v_{x}
\label{(2.2)}
\end{equation}
is the Jacobian of the map $(x,y) \rightarrow (u,v)$, $\Phi(x,y)$ is the upper limit
\begin{equation}
\Phi(x,y) \equiv {\overline {\lim_{\bigtriangleup z \to 0}}}
\left | {\bigtriangleup w \over \bigtriangleup z} \right | ,
\label{(2.3)}
\end{equation}
$\varphi(x,y)$ is the lower limit
\begin{equation}
\varphi(x,y) \equiv {\underline {\lim_{\bigtriangleup z \to 0}}}
\left | {\bigtriangleup w \over \bigtriangleup z} \right | ,
\label{(2.4)}
\end{equation}
then $J \geq 0$ almost everywhere in $A$, and there exists a positive real number 
$\mu \in ]0,1]$ such that
\begin{equation}
\varphi(x,y) \geq \mu \Phi(x,y)
\label{(2.5)}
\end{equation}
almost everywhere in $A$. The function $w$ is here said to be pseudoholomorphic\footnote{The work
in Refs. \cite{Caccioppoli1952,Caccioppoli1953} uses actually the nomenclature {\it pseudoanalytic},
but we prefer to speak of pseudoholomorphic functions, to avoid confusion with the local theory of
pseudoanalytic functions, which relies instead upon generalized Cauchy-Riemann equations \cite{Courant1962}}. 
of parameter $\mu$ \cite{Caccioppoli1952,Caccioppoli1953} 
Every value of $\mu \leq 1$ corresponds to a class $C_{\mu}$ of pseudholomorphic functions; in
particular, if $\mu=1$, $C_{1}$ is the class of holomorphic functions.

Outside of the holomorphic framework, the increment ratio 
${\bigtriangleup w \over \bigtriangleup z}$ has indeed a rich structure because, upon defining
\begin{equation}
m \equiv {\bigtriangleup y \over \bigtriangleup x}
\label{(2.6)}
\end{equation}
one has
\begin{equation}
{\bigtriangleup w \over \bigtriangleup z}
={u_{x}+{\rm i}v_{x} +m(u_{y}+{\rm i}v_{y}) \over (1+{\rm i}m)}
+{(\varepsilon_{1}+{\rm i}\varepsilon_{2})+m(\varepsilon_{3}+{\rm i}\varepsilon_{4})
\over (1+{\rm i}m)},
\label{(2.7)}
\end{equation}
where $\varepsilon_{i}$ tends to zero as $\bigtriangleup z \rightarrow 0$, for all
$i=1,2,3,4$. Multiplication and division by $(1-{\rm i}m)$ on the right-hand side yields therefore
\begin{equation}
{\bigtriangleup w \over \bigtriangleup z}=A+{\rm i}B+{\rm O}(\varepsilon),
\label{(2.8)}
\end{equation}
where 
\begin{equation}
A \equiv {(u_{x}+m u_{y})+m(v_{x}+m v_{y}) \over (1+m^{2})},
\label{(2.9)}
\end{equation}
\begin{equation}
B \equiv {(v_{x}+m v_{y})-m(u_{x}+m u_{y}) \over (1+m^{2})},
\label{(2.10)}
\end{equation}
and hence
\begin{equation}
\lim_{\bigtriangleup z \to 0} \left | {\bigtriangleup w \over \bigtriangleup z} \right |^{2}
=A^{2}+B^{2}.
\label{(2.11)}
\end{equation}
Now a patient calculation shows that, in the expression $(1+m^{2})^{2}(A^{2}+B^{2})$, a
cancellation occurs and some reassembling can be made, so that eventually
\cite{Hedrick1933}
\begin{equation}
A^{2}+B^{2}=(1+m^{2})^{-1}(E+2Fm+Gm^{2}) \equiv r(m),
\label{(2.12)}
\end{equation}
where, according to a standard notation, we have set
\begin{equation}
E \equiv (u_{x})^{2}+(v_{x})^{2}, \;
G \equiv (u_{y})^{2}+(v_{y})^{2}, \;
F \equiv u_{x}u_{y}+v_{x}v_{y}.
\label{(2.13)}
\end{equation}
To study the maximum or minimum of $A^{2}+B^{2}$ as a function of 
the real variable $m$, we have to evaluate its
first derivative, which vanishes if $m$ solves the algebraic equation of second degree
\begin{equation}
F m^{2}+(E-G)m-F=0,
\label{(2.14)}
\end{equation}
solved by
\begin{equation}
m=-{(E-G) \over 2F} \pm {1 \over 2F} \sqrt{(E-G)^{2}+4F^{2}}=m_{\pm}.
\label{(2.15)}
\end{equation}

It is clear from (\ref{(2.15)}) that for any generic E, F and G, $m_{+}$ is positive whereas $m_{-}$ is negative.
Once we take second derivative of  the function $A^{2}+B^{2}$ with respect to $m$, after a little algebra,
we obtain the simple expression 
\begin{equation}
\frac{\partial^{2}}{\partial m^{2}} (A^{2}+B^{2}) = - \frac{2F}{m(1+m^{2})}
\label{(2.16)}
\end{equation}
and hence it is clear that $\frac{\partial^{2}}{\partial m^{2}} (A^{2}+B^{2})|_{m_{+}} < 0$ whereas
$\frac{\partial^{2}}{\partial m^{2}} (A^{2}+B^{2})|_{m_{-}} > 0$.

Interestingly, setting for convenience
\begin{equation}
\omega \equiv \sqrt{(E-G)^{2}+4F^{2}}=\sqrt{(E+G)^{2}-4J^{2}},
\label{(2.17)}
\end{equation}
one finds from (2.12) that \cite{Hedrick1933}
\begin{equation}
r_{+}=r(m_{+})={E+ 2 F m_{+} +G m_{+}^{2} \over (1+m_{+}^{2})}
={(E+G) \over 2}+{\omega \over 2},
\label{(2.18)}
\end{equation}
\begin{equation}
r_{-}=r(m_{-})={E+ 2 F m_{-} +G m_{-}^{2} \over (1+m_{-}^{2})}
={(E+G) \over 2}-{\omega \over 2}.
\label{(2.19)}
\end{equation}
In other words, the maximum and minimum values of $r$ are themselves solutions of the
quadratic equation \cite{Hedrick1933}
\begin{equation}
\rho^{2}-(E+G)\rho+J^{2}=0,
\label{(2.20)}
\end{equation}
and the upper and lower limit (2.3) and (2.4) turn out to obey the identities 
\cite{Caccioppoli1952,Caccioppoli1953}
\begin{equation}
2 \Phi^{2}=E+G+ \omega, \; 
2 \varphi^{2}=E+G- \omega,
\label{(2.21)}
\end{equation}
\begin{equation}
\Phi^{2}+\varphi^{2}=E+G, \;
\Phi^{2} \varphi^{2}={1 \over 4}[(E+G)^{2}-\omega^{2}]=J^{2}
=EG-F^{2},
\label{(2.22)}
\end{equation}
jointly with the inequalities
\begin{equation}
\mu J \leq \varphi^{2} \leq \Phi^{2} \leq {J \over \mu}, \;
J \geq {\mu \over (1+\mu^{2})}(\Phi^{2}+\varphi^{2}).
\label{(2.23)}
\end{equation}
By virtue of assumption (i), one can build a sequence $(u_{n}(x,y),v_{n}(x,y))$ of pairs of functions,
describing flat surfaces, such that \cite{Caccioppoli1952}
\begin{equation}
\lim_{n \to \infty}u_{n}(x,y)=u(x,y), \;
\lim_{n \to \infty}v_{n}(x,y)=v(x,y),
\label{(2.24)}
\end{equation}
uniformly in every closed portion of the field $A$, and such that their partial derivatives with respect
to $x$ and $y$ converge in mean of order $2$; $u_{n}$ and $v_{n}$ being functions as smooth as desired,
or even polynomials. Hence it follows that on the surface $S$ associated to $u(x,y)$ and $v(x,y)$, areas 
and lengths have the classical expressions. This means that, if $D$ and $L$ are a domain and a line within
$A$, respectively, the area $\tau D$ and the length $\tau L$ read as
\begin{equation}
\int \int_{D}J \; {\rm d}x \; {\rm d}y
=\int \int_{D}\sqrt{EG-F^{2}} \; {\rm d}x \; {\rm d}y,
\label{(2.25)}
\end{equation}
\begin{equation}
\tau L = \int_{L}\sqrt{E {\rm d}x^{2}+2F {\rm d}x \; {\rm d}y+G {\rm d}y^{2}}.
\label{(2.26)}
\end{equation}
With this nomenclature, $\tau$ is the plane transformation of $z$ into $w$ defined by the equations
$$
\tau: \; u=u(x,y), \; v=v(x,y)
$$
which in turn describe a flat surface $S$ carried by the plane $w$. Such a map $\tau$ is said to be
a pseudo-conformal transformation with parameter $\mu$, i.e. a pseudo-conformal representation
of $S$ upon $A$. To every point $z_{0}$ of the field $A$ there corresponds a point $P(z_{0})$ of $S$,
having as trace on the plane $w$ the point $w_{0}=w(z_{0})$.

\subsection{Quasi-conformal maps for real and complex manifolds}

For the case of real manifolds we recall, following Ref. \cite{Donaldson1989}, that for any
pseudo-group of homeomorphisms of Euclidean space one can define the corresponding category of manifolds.
Thus, the full pseudo-group of homeomorphisms, the subgroup of smooth diffeomorphisms, the pseudo-group
of quasiconformal maps and the pseudo-group of Lipschitz maps give rise to the theory of topological
manifolds, $C^{\infty}$ manifolds, quasi-conformal manifolds and Lipschitz manifolds, respectively. 
In particular, a homeomorphism $\varphi: D \rightarrow {\bf R}^{n}$ is $K$-quasiconformal if,
for all $x \in D$,
$$
\lim_{r \to 0} {\rm sup} {{\rm max} |\varphi(y)-\varphi(x)| \over {\rm min}
|\varphi(y)-\varphi(x)|} \leq K,
$$
with
$$
|y-x|=\sqrt{\sum_{k=1}^{n}(y_{k}-x_{k})^{2}}=r.
$$
The map $\varphi$ is quasiconformal if it is $K$-quasiconformal for some $K \geq 1$. This range of
values of $K$ is responsible for distortion of relative distances of nearby points by a bounded factor.

On the other hand, for functions of complex variable, 
the concept of quasiconformal maps \cite{Bers1957} was introduced by
Gr\"{o}tzsch, who considered homeomorphisms (2.1) with a positive Jacobian (2.2). Such a map
takes infinitesimal circles (cf. the above remarks on distortion of relative distances) 
into infinitesimal ellipses, and is called {\it quasi-conformal}
if the eccentricity of these ellipses is uniformly bounded. This condition can be expressed
analytically by either of the three equivalent differential inequalities in (2.23)
(our $\mu$ parameter is the inverse of the $Q$ parameter used by Bers). This property is
conformally invariant: if $w=w(z)$ has it, so does the function $U(\zeta)=F \left \{
w [f(\zeta)] \right \}$, where $F$ and $f$ are conformal mappings.

Bers \cite{Bers1957} calls a function $w(z)$ as in (2.1) quasi-conformal if it is of the form
\begin{equation}
w(z)=f[\chi(z)],
\label{(2.27)}
\end{equation}
where $\chi$ is a quasi-conformal homeomorphism and $f$ is an holomorphic function. This
definition is suggested by an important result of Morrey \cite{Morrey1938} according to which,
if $w(z)$ is a quasi-conformal function defined in the unit disk, then it admits the
representation (2.27) where $\zeta=\chi(z)$ is a homeomorphism of the set $|z| \leq 1$
onto the set $|\zeta| \leq 1$ with $\chi(0)=0,\chi(1)=1$, which satisfies together with
its inverse $\chi^{-1}$ a uniform Holder condition, and where $f(\zeta)$ is a holomorphic
function of the complex variable $\zeta$, $|\zeta|<1$.

It is a profound result of the work in Ref. \cite{Bers1957}, which was inspired also by the
work of Mori \cite{Mori1957}, that the geometric definition relying upon Eq. (2.27) is equivalent
to the analytic definition given at the beginning of this section, inspired by Morrey,
Caccioppoli, Bers and Nirenberg, according to which a continuous function 
$w(z)=u(x,y)+{\rm i}v(x,y)$ in a domain $D$ is quasi-conformal if it has $L_{2}$ derivatives
satisfying the differential inequalities (2.23) almost everywhere. By equivalent we here
mean that the geometric and analytic definition imply each other \cite{Bers1957}.
Note that, in light of the results in Refs. \cite{Caccioppoli1952,Caccioppoli1953,Bers1957},
our pseudoholomorphic functions are also quasi-conformal maps.

\section{Pseudoholomorphic functions of $n$ complex variables}

Let us proceed now by generalizing the results of previous section to $n$ complex variables.

For $n$ complex variables $z_{1}$, $z_{2}$,........$z_{n}$ with  
$z_{k}=x_{k}+{\rm i}y_{k}$, $\forall k = 1,2,...,n$  $x_{k} \in {\bf R}, y_{k} \in
{\bf R}$, let $w$ be the continuous function with image
\begin{equation}
w(z_{1}, z_{2},....,z_{n})=u(x_{1}, x_{2},....,x_{n},y_{1}, y_{2},....,y_{n})
+{\rm i}v(x_{1}, x_{2},....,x_{n},y_{1}, y_{2},.....,y_{n}).
\label{(3.1)}
\end{equation}
Let the functions $u(\left\{x_{k}\right\}, \left\{y_{k}\right\}),v(\left\{x_{k}\right\}, 
\left\{y_{k}\right\})$ satisfy the following assumptions:
\vskip 0.3cm
\noindent
(i) $u(\left\{x_{k}\right\}, \left\{y_{k}\right\})$ and $v(\left\{x_{k}\right\}, \left\{y_{k}\right\})$ 
are absolutely continuous in $\left\{x_{k}\right\}$ and $\left\{y_{k}\right\}$ for almost all values of  
$\left\{y_{k}\right\}, \left\{x_{k}\right\}$ respectively, while their first derivatives 
$\left\{u_{x_k}\right\}, \left\{u_{y_k}\right\},\left\{v_{x_k}\right\},\left\{v_{y_k}\right\}$ 
are square-integrable in every internal portion of the domain of definition.
\vskip 0.3cm
\noindent
(ii) If
\begin{equation}
J_k={\partial (u,v) \over \partial (x_k,y_k)}=u_{x_k}v_{y_k}-u_{y_k}v_{x_k},
\label{(3.2)}
\end{equation}
$\Phi_k$ is the upper limit
\begin{equation}
\Phi_k \equiv {\overline {\lim_{\bigtriangleup z_k \to 0}}}
\left | {\bigtriangleup w \over \bigtriangleup z_k} \right | ,
\label{(3.3)}
\end{equation}
$\varphi_k$ is the lower limit
\begin{equation}
\varphi_k \equiv {\underline {\lim_{\bigtriangleup z_k \to 0}}}
\left | {\bigtriangleup w \over \bigtriangleup z_k} \right | ,
\label{(3.4)}
\end{equation}
then $J_k \geq 0$ almost everywhere in $A$, and there exists a positive real number 
$\mu_k \in  ]0,1]$ such that
\begin{equation}
\varphi_k \geq \mu_k \Phi_k
\label{(3.5)}
\end{equation}
almost everywhere in $A$. We then say that the function $w$ is  
pseudoholomorphic of parameter $\mu_k$ for all $k=1,2,...,n$.

Outside of the holomorphic framework, the increment ratio 
${\bigtriangleup w \over \bigtriangleup z_k}$ has indeed a rich structure 
as we know already from Sect. 2 because, upon defining
\begin{equation}
m_k \equiv {\bigtriangleup y_k \over \bigtriangleup x_k}
\label{(3.6)}
\end{equation}
one has $n$ increment ratios
\begin{equation}
{\bigtriangleup w \over \bigtriangleup z_k}
={u_{x_k}+{\rm i}v_{x_k} +m_k(u_{y_k}+{\rm i}v_{y_k}) \over (1+{\rm i}m_k)}
+{(\varepsilon_{1}+{\rm i}\varepsilon_{2})_{k}+m_k(\varepsilon_{3}+{\rm i}\varepsilon_{4})_{k}
\over (1+{\rm i}m_k)},
\label{(3.7)}
\end{equation}
where $(\varepsilon_{i})_{k}$ tends to zero as $\bigtriangleup z_k \rightarrow 0$, for all
$i=1,2,3,4$ and for all $k=1,2,.....,n$. Multiplication and division by 
$(1-{\rm i}m_k)$ on the right-hand side of (3.7) yields therefore
\begin{equation}
{\bigtriangleup w \over \bigtriangleup z_k}=A_k+{\rm i}B_k+{\rm O}(\varepsilon),
\label{(3.8)}
\end{equation}
where 
\begin{equation}
A_k \equiv {(u_{x_k}+m_k u_{y_k})+m_k(v_{x_k}+m_k v_{y_k}) \over (1+m_k^{2})},
\label{(3.9)}
\end{equation}
\begin{equation}
B_k \equiv {(v_{x_k}+m_k v_{y_k})-m_k(u_{x_k}+m_k u_{y_k}) \over (1+m_k^{2})},
\label{(3.10)}
\end{equation}
and hence
\begin{equation}
\lim_{\bigtriangleup z_k \to 0} \left | {\bigtriangleup w \over \bigtriangleup z_k} \right |^{2}
=A_k^{2}+B_k^{2}.
\label{(3.11)}
\end{equation}
With the help of the same calculations leading to Eq. (2.12) we now find
\begin{equation}
A_k^{2}+B_k^{2}=(1+m_k^{2})^{-1}(E_k+2F_km_k+G_km_k^{2}) \equiv r(m_{k}),
\label{(3.12)}
\end{equation}
where, according to our own notation, we have set
\begin{equation}
E_k \equiv (u_{x_k})^{2}+(v_{x_k})^{2}, \;
G_k \equiv (u_{y_k})^{2}+(v_{y_k})^{2}, \;
F_k \equiv u_{x_k}u_{y_k}+v_{x_k}v_{y_k}.
\label{(3.13)}
\end{equation}
To study the maximum or minimum of $A_k^{2}+B_k^{2}$ as a function of $m_k$, we have to evaluate its
first derivative, which vanishes if $m_k$ solves the algebraic equation of second degree
\begin{equation}
F_k m_k^{2}+(E_k-G_k)m_k-F_k=0,
\label{3.14)}
\end{equation}
solved by
\begin{equation}
m_k=-{(E_k-G_k) \over 2F_k} \pm {1 \over 2F_k} \sqrt{(E_k-G_k)^{2}+4F_k^{2}}=(m_{\pm})_{k}.
\label{(3.15)}
\end{equation}
It is clear from (\ref{(3.15)}) that for any generic $E_k$, $F_k$ and $G_k$, $(m_{+})_{k}$ 
is positive whereas $(m_{-})_{k}$ is negative.
And then once we take second derivative of  the function $A_k^{2}+B_k^{2}$ with respect to $m_k$, 
we find, as in Eq. (2.16),
\begin{equation}
\frac{\partial^{2}}{\partial m_k^{2}} (A_k^{2}+B_k^{2}) = - \frac{2F_k}{m_k(1+m_k^{2})},
\label{(3.16)}
\end{equation}
and hence it is clear that $\frac{\partial^{2}}{\partial m_k^{2}} (A_k^{2}+B_k^{2})|_{(m_{+})_{k}} < 0$ and
$\frac{\partial^{2}}{\partial m^{2}} (A^{2}+B^{2})|_{(m_{-})_{k}} > 0$.

Interestingly, setting for convenience
\begin{equation}
\omega_k \equiv \sqrt{(E_k-G_k)^{2}+4F_k^{2}}=\sqrt{(E_k+G_k)^{2}-4J_k^{2}},
\label{(3.17)}
\end{equation}
one finds from (3.12) that (cf. (2.18) and (2.19))
\begin{equation}
(r_{+})_{k}=r((m_{+})_{k})={E_k+ 2 F_k (m_{+})_{k} +G_k (m_{+}^{2})_{k} \over (1+(m_{+}^{2})_{k})}
={(E_k+G_k) \over 2}+{\omega_k \over 2},
\label{(3.18)}
\end{equation}
\begin{equation}
(r_{-})_{k}=r((m_{-})_{k})={E_k+ 2 F_k (m_{-})_{k} +G_k (m_{-}^{2})_{k} \over (1+(m_{-}^{2}) _{k})}
={(E_k+G_k) \over 2}-{\omega_k \over 2}.
\label{(3.19)}
\end{equation}
In other words, the maximum and minimum values of $r_k$ are themselves solutions of the
quadratic equation (cf. Eq. (2.20))
\begin{equation}
\rho^{2}-(E_k+G_k)\rho+J_k^{2}=0,
\label{(3.20)}
\end{equation}
and the upper and lower limit (3.3) and (3.4) turn out to obey the relations
(cf. (2.21)-(2.23))
\begin{equation}
2 \Phi_k^{2}=E_k+G_k+ \omega_k, \; 
2 \varphi_k^{2}=E_k+G_k- \omega_k,
\label{(3.21)}
\end{equation}
\begin{equation}
\Phi_k^{2}+\varphi_k^{2}=E_k+G_k, \;
\Phi_k^{2} \varphi_k^{2}={1 \over 4}[(E_k+G_k)^{2}-\omega_k^{2}]=J_k^{2}
=E_kG_k-F_k^{2},
\label{(3.22)}
\end{equation}
\begin{equation}
\mu_k J_k \leq \varphi_k^{2} \leq \Phi_k^{2} \leq {J_k \over \mu_k}, \;
J_k \geq {\mu_k \over (1+\mu_k^{2})}(\Phi_k^{2}+\varphi_k^{2}).
\label{(3.23)}
\end{equation}
As one can see, for the $n$-variable case the functions $u$ and $v$ are severely constrained, since there is
an $n$-tuple of conditions to be satisfied. 

\section{Discwise quasi-conformal maps with $n$ complex variables}

In the attempt of providing examples of the functions fulfilling our definition in Sect. 3, we here
present a review of part of the work in Ref. \cite{Hitotumatu1959}, devoted to the investigation
of quasi-conformal functions of several complex variables.

A function $f$ of $n$ complex variables $z_{1},...,z_{n}$ defined in a domain $D$ is said to be
a {\it discwise-quasi-conformal function} with dilatation $K={1 \over \mu}$ if
\vskip 0.3cm
(i) $f$ is of class $C^{1}$ in $D$, and
\vskip 0.3cm
(ii) $f$ is $K$-quasi-conformal (i.e. quasi-conformal with dilatation
parameter $K={1 \over \mu}$) on each holomorphic plane.
\vskip 0.3cm
\noindent
The meaning of condition (ii) is as follows. Once we have a linear map 
\begin{equation}
z_{1}=a_{1}t+b_{1},...,z_{n}=a_{n}t+b_{n},
\label{(4.1)}
\end{equation}
the $a$'s and $b$'s being complex constants while $t$ is a complex variable,
defined on the unit disc $\left \{ |t| < 1 \right \}$ and whose image lies completely
in $D$, the composite function 
\begin{equation}
f(t)=f(a_{1}t+b_{1},...,a_{n}t+b_{n})
\label{(4.2)}
\end{equation}
is always a $K$-quasi-conformal function in the unit disc $\left \{ | t | < 1
\right \}$. The function $f$ is then said to be $K$-{\it discwise-quasi-conformal},
following Ref. \cite{Hitotumatu1959}. In particular, a $1$-discwise-quasi-conformal function 
is a holomorphic function in $D$. A number of important properties are found to hold, and
they are as follows \cite{Hitotumatu1959}.
\vskip 0.3cm
\noindent
{\bf Theorem 4.1.} If $f$ is discwise quasi-conformal in a bounded domain $D$ and continuous
also on the boundary $\partial D$ of $D$, then the maximum principle holds, according to which
$$
{\rm sup} \left \{|f|;D \right \}
= {\rm sup} \left \{ |f| ; \partial D \right \}.
$$
\vskip 0.3cm
\noindent
{\bf Theorem 4.2.} If, on the unit disc $\left \{ |t| < 1 \right \}$, one considers the 
holomorphic map 
\begin{equation}
z_{1}=\varphi_{1}(t),...,z_{n}=\varphi_{n}(t),
\label{(4.3)}
\end{equation}
whose image lies completely in $D$, the composite function
\begin{equation}
{\hat f}(t)=f(\varphi_{1}(t),...,\varphi_{n}(t))
\label{(4.4)}
\end{equation}
is a $K$-quasi-conformal function in the unit disc. As a corollary, the concept of
$K$-discwise-quasi-conformal function is invariant under holomorphic transformations. This means
that, if $f(z_{1},...,z_{n})$ is $K$-discwise-quasi-conformal in $D$, and if
\begin{equation}
z_{1}=\varphi_{1}(w_{1},...,w_{m}),...,
z_{n}=\varphi_{n}(w_{1},...,w_{m})
\label{(4.5)}
\end{equation}
is a holomorphic transformation from a domain $B$ in
$(w_{1},...,w_{m})$-space into $D$, then the composite function
\begin{equation}
F(w_{1},...,w_{m})=f(\varphi_{1}(w_{1},...,w_{m}),...,\varphi_{n}(w_{1},...,w_{m}))
\label{(4.6)}
\end{equation}
is again $K$-discwise-quasi-conformal in $B$.
One can therefore define a $K$-discwise-quasi-conformal function as a smooth function which is
$K$-quasi-conformal on every holomorphic surface. This holds not only in a domain $D$, but also
for an arbitrary set, in particular on an analytic subset in the $(z_{1},...,z_{n})$-space.

Another peculiar property is that the sum of two $K$-discwise-quasi-conformal functions is not
always $K$-discwise-quasi-conformal. For example, each of the functions 
$$
2z_{1}+2z_{2}+{\overline z}_{1}+{\overline z}_{2} \; {\rm and}
-z_{1}-z_{2},
$$
is $K$-discwise-quasi-conformal, but their sum, being equal to
$$
2({\rm Re}z_{1}+{\rm Re}z_{2}),
$$
is not $K$-discwise-quasi-conformal, because quasi-conformal functions cannot take only real
values (having to provide an open mapping) unless they are a constant. 
\vskip 0.3cm
\noindent
{\bf Theorem 4.3.} At every ordinary point, where at least one of the first partial derivatives with 
respect to $z_{1},...,z_{n}$ does not vanish, the real and imaginary parts of a
$K$-discwise-quasi-conformal function $f=u+{\rm i}v$ satisfy the following system of
partial differential equations identically:
\begin{equation}
{\partial (u,v) \over \partial (x_{j},x_{k})}
={\partial (u,v) \over \partial (y_{j},y_{k})}, \;
{\partial (u,v) \over \partial (x_{j},y_{k})}
={\partial (u,v) \over \partial (x_{k},y_{j})}, \;
j,k=1,2,...,n,
\label{(4.7)}
\end{equation}
where the notation for the independent variables is the same as in (3.1).
\vskip 0.3cm
\noindent
{\bf Theorem 4.4.} If $f$ is a $K$-discwise-quasi-conformal function of two complex variables
$z_{1},z_{2}$ in a domain $D$, then at every ordinary point $(z_{1}^{0},z_{2}^{0})$ the set
given by the equation 
\begin{equation}
f(z_{1},z_{2})=f(z_{1}^{0},z_{2}^{0})={\rm constant}
\label{(4.8)}
\end{equation}
is a two-dimensional holomorphic surface.
\vskip 0.3cm
\noindent
{\bf Theorem 4.5.} If $\kappa$ is a given continuous function of $2n$ real variables 
$x_{1},y_{1},...,x_{n},y_{n}$ in a domain $D$, and if the modulus of $\kappa$ is bounded
by a constant $k_{0}<1$, then a solution of class $C^{1}$ in $D$ of the system of partial 
differential equations (4.7) yields a $K$-discwise-quasi-conformal function
$f=u+{\rm i}v$ in $D$, whose dilatation $K$ is given by
\begin{equation}
K={(1 + k_{0}) \over (1-k_{0})}, \;
|\kappa| \leq k_{0}<1.
\label{(4.9)}
\end{equation}

The results here recalled are helpful in understanding properties and limits of the
quasi-conformal  pseudoholomorphic framework that we are investigating, and can be
compared with the different perspectives considered in Refs.
\cite{Toki1954,Storvick1957,Koohara1971,Fryant1981,Tutschke2007}. 

\section{The real manifold associated with a pseudoholomorphic function}

Since the conditions (2.23) are differential inequalities, it is at first sight problematic 
to define a pseudoholomorphic manifold, even just in the case of one complex dimension.
However, we may point out that, if $\varphi^{2} > \mu J$ is fulfilled (we rule out the case of
equalities, for which we refer the reader to the Appendix), then also $\Phi^{2} > \mu J$ is
fulfilled, because $\Phi^{2}> \varphi^{2}$. Moreover, for the same reason, if 
$\Phi^{2}< {J \over \mu}$ holds, {\it a fortiori} $\varphi^{2}< {J \over \mu}$ holds as well.
We can therefore consider three positive-definite functions $\alpha,\beta$ and $\gamma$
depending on $(\mu;x,y)$, such that (2.23) is re-expressed in the form  
\begin{equation}
\varphi^{2}-\mu J =\alpha > 0,
\label{(5.1)}
\end{equation}
\begin{equation}
{J \over \mu}-\Phi^{2}=\beta >0,
\label{(5.2)}
\end{equation}
\begin{equation}
J-{\mu \over (1+\mu^{2})}(\Phi^{2}+\varphi^{2}) =\gamma > 0.
\label{(5.3)}
\end{equation}
Note also that (5.2) leads to $\mu \Phi^{2}=J-\mu \beta$, while (5.1) implies
that $\mu \varphi^{2}=\mu^{2}J+\mu \alpha$. Thus, on the left-hand side of (5.3),
we obtain eventually exact cancellation of terms proportional to $J$, so that
\begin{equation}
-{\mu(\alpha-\beta)\over (1+\mu^{2})}=\gamma >0 \Longrightarrow 
\alpha(\mu;x,y) < \beta(\mu;x,y).
\label{(5.4)}
\end{equation}
We can now exploit (2.17) and (2.21) to re-express (5.1) and (5.2) in the form
\begin{equation}
E+G-\sqrt{(E+G)^{2}-4 J^{2}}-2 \mu J -2 \alpha(\mu;x,y)=0,
\label{(5.5)}
\end{equation}
\begin{equation}
2{J \over \mu}-E-G-\sqrt{(E+G)^{2}-4J^{2}}-2 \beta(\mu;x,y)=0.
\label{(5.6)}
\end{equation}
Within this framework, the desired pseudoholomorphic manifold is defined {\it implicitly}
by the nonlinear equations (5.5) and (5.6) (actually by a single equation equivalent to
them, see below), a procedure which is often used in simpler cases
\cite{Rheinboldt1996} in the literature.  Upon using (2.2) and (2.13), we here look for 
$u(x,y),v(x,y),\alpha(\mu,x,y)$ and $\beta(\mu,x,y)$ such that (5.5) and (5.6) hold,
and the consistency condition (5.4) is fulfilled. At this stage, we can express
$\sqrt{(E+G)^{2}-4J^{2}}$ from Eq. (5.5), and then insert the result into Eq. (5.6). This
leads to a nonlinear equation without square root, reading as
\begin{equation}
{(\mu^{2}+1)\over \mu}J-E-G+\beta -\alpha=0,
\label{(5.7)}
\end{equation}
which is nothing but Eq. (5.3). In particular, we may look for $u(x,y)$ and $v(x,y)$ in the form
\begin{equation}
u(x,y)=Ax^{2}+Bxy+Cy^{2},
\label{(5.8)}
\end{equation}
\begin{equation}
v(x,y)={\widetilde A}x^{2}+{\widetilde B}xy+{\widetilde C}y^{2}.
\label{(5.9)}
\end{equation}
This ansatz turns Eq. (5.7) into the form
\begin{equation}
A_{\mu}x^{2}+B_{\mu}xy+C_{\mu}y^{2}=(\alpha-\beta)<0,
\label{(5.10)}
\end{equation}
where the coefficients are given by
\begin{equation}
A_{\mu} \equiv 2{(\mu^{2}+1)\over \mu}(A{\widetilde B}-B{\widetilde A})
-(4A^{2}+4{\widetilde A}^{2}+B^{2}+{\widetilde B}^{2}),
\label{(5.11)}
\end{equation}
\begin{equation}
B_{\mu} \equiv 4 \left[{(\mu^{2}+1)\over \mu}(A{\widetilde C}-C{\widetilde A})-(A+C)B
-({\widetilde A}+{\widetilde C}){\widetilde B}\right],
\label{(5.12)}
\end{equation}
\begin{equation}
C_{\mu} \equiv 2{(\mu^{2}+1)\over \mu}(B{\widetilde C}-C{\widetilde B})
-(4C^{2}+4{\widetilde C}^{2}+B^{2}+{\widetilde B}^{2}).
\label{(5.13)}
\end{equation}

For simplicity, we may study first the particular case when (5.10) is fulfilled with
$B_{\mu}=0,A_{\mu}<0,C_{\mu}<0$. Indeed, $B_{\mu}$ vanishes if $A=-C$ and 
${\widetilde A}=-{\widetilde C}$, and in such a case $A_{\mu}=C_{\mu}$. The quadratic form
in (5.10) is then negative-definite if $A,{\widetilde A},B,{\widetilde B}$ are chosen so as
to satisfy the condition
\begin{equation}
A{\widetilde B}-B{\widetilde A} < {\mu \over 2(\mu^{2}+1)}
\Bigr(4A^{2}+4{\widetilde A}^{2}+B^{2}+{\widetilde B}^{2}\Bigr).
\label{(5.14)}
\end{equation}
In the general case, the condition (5.10) is fulfilled if the quadratic form
$-A_{\mu}x^{2}-B_{\mu}xy-C_{\mu}y^{2}$ is positive-definite. For this purpose, we have to study
the eigenvalues $\lambda$ which are the roots of the following equation of second degree:
\begin{equation}
{\rm det} \left(\begin{matrix} (-A_{\mu}-\lambda) & -{1 \over 2}B_{\mu} \cr
-{1 \over 2}B_{\mu} & (-C_{\mu}-\lambda) \end{matrix}\right)=0.
\label{(5.15)}
\end{equation}
This equation is solved by the two roots
\begin{equation}
\lambda=-{1 \over 2}(A_{\mu}+C_{\mu}) \pm {1 \over 2} \sqrt{\delta_{ABC}}
=\lambda_{1,2},
\label{(5.16)}
\end{equation}
where
\begin{equation}
\delta_{ABC}=(A_{\mu}+C_{\mu})^{2}-4 \left(A_{\mu}C_{\mu}-{1 \over 4}B_{\mu}^{2}\right)
=(A_{\mu}-C_{\mu})^{2}+B_{\mu}^{2} > 0,
\label{(5.17)}
\end{equation}
and $\lambda_{1}$ (resp. $\lambda_{2}$) corresponds to the $+$ (resp. $-$) sign in front of
$\sqrt{\delta_{ABC}}$. Since $\delta_{ABC}$ is a sum of squares of real numbers, we are
guaranteed that both roots are real. Moreover, 
\begin{equation}
\lambda_{1}-\lambda_{2}=\sqrt{\delta_{ABC}} > 0,
\label{(5.18)}
\end{equation}
and hence positivity of $\lambda_{2}$ ensures positivity of $\lambda_{1}$ as well. By virtue
of (5.16), positivity of $\lambda_{2}$ implies that
\begin{equation}
-\sqrt{\delta_{ABC}} > (A_{\mu}+C_{\mu}) \Longrightarrow (A_{\mu}+C_{\mu}) <0
\Longrightarrow \sqrt{\delta_{ABC}} < |A_{\mu}+C_{\mu}|.
\label{(5.19)}
\end{equation}

In other words, (5.10) is always fulfilled provided that the following two conditions hold:
\begin{equation}
(A_{\mu}+C_{\mu}) < 0,
\label{(5.20)}
\end{equation}
\begin{equation}
\sqrt{(A_{\mu}-C_{\mu})^{2}+B_{\mu}^{2}} < |A_{\mu}+C_{\mu}|,
\label{(5.21)}
\end{equation}
where $A_{\mu},B_{\mu},C_{\mu}$ are displayed in (5.11)-(5.13).

We have still to understand whether our Eq. (5.10) can be seen to define a manifold. For this purpose, 
we recall a basic theorem \cite{Schwartz1967}, according to which, for a set $V$ of an affine space
$E$ to be an hypersurface (i.e. a manifold of dimension $N-1$) of class $C^{m}$ of $E$, it is
necessary and sufficient that, $\forall a \in V$, there exists a neighbourhood $U(a)$ of $a$ within
$E$ and a scalar function $H$, defined in $U$ and of class $C^{m}$, such that, in a reference
system, the partial derivatives ${\partial H \over \partial x^{i}}(a), i=1,2,...,N$ are not all
simultaneously vanishing. In our case, the scalar function pertaining to Eq. (5.10) is
the $1$-parameter family
\begin{equation}
H(\mu; x,y) \equiv A_{\mu}x^{2}+B_{\mu}xy+C_{\mu}y^{2}
+\beta(\mu;x,y)-\alpha(\mu;x,y),
\label{(5.22)}
\end{equation}
with gradient having components
\begin{equation}
{\partial H \over \partial x}=2 A_{\mu}x+B_{\mu}y+{\partial \over \partial x}(\beta-\alpha),
\label{(5.23)}
\end{equation}
\begin{equation}
{\partial H \over \partial y}=B_{\mu}x+2C_{\mu}y
+{\partial \over \partial y}(\beta-\alpha).
\label{(5.24)}
\end{equation}
We may therefore obtain a one-dimensional real manifold, i.e. a curve $\Gamma$ defined
implicitly as the set
\begin{equation}
\Gamma \equiv \left \{ a=(x,y): ({\rm grad} \; H)(a) \not =0 \right \}.
\label{(5.25)}
\end{equation}
The gradient of $(\beta-\alpha)$ plays a crucial role in ensuring that the gradient of $H$
does not vanish (otherwise we would have the singular point $(x=0,y=0)$, and hence
no manifold exits). 

Since Eq. (5.10) is just a particular case of Eq. (5.7), we should also consider the general form of
Eq. (5.22), i.e.
\begin{equation}
H(\mu;x,y) \equiv 
(u_{x})^{2}+(v_{x})^{2}+(u_{y})^{2}+(v_{y})^{2}
-{(\mu^{2}+1)\over \mu}(u_{x}v_{y}-u_{y}v_{x})+\alpha-\beta.
\label{(5.26)}
\end{equation}
We are therefore looking for a $1$-parameter family of curves defined as in (5.25), where
the gradient of $H$ has components
\begin{eqnarray}
{\partial H \over \partial x}&=& 
\left(2u_{x}-{(\mu^{2}+1)\over \mu}v_{y}\right)u_{xx}
+\left(2v_{x}+{(\mu^{2}+1)\over \mu}u_{y}\right)v_{xx}
\nonumber \\
&+& \left(2u_{y}+{(\mu^{2}+1)\over \mu}v_{x}\right)u_{xy}
+\left(2v_{y}-{(\mu^{2}+1)\over \mu}u_{x}\right)v_{xy}
\nonumber \\
&+& \alpha_{x}-\beta_{x},
\label{(5.27)}
\end{eqnarray}
\begin{eqnarray}
{\partial H \over \partial y}&=& 
\left(2u_{y}+{(\mu^{2}+1)\over \mu}v_{x}\right)u_{yy}
+\left(2v_{y}-{(\mu^{2}+1)\over \mu}u_{x}\right)v_{yy}
\nonumber \\
&+& \left(2u_{x}-{(\mu^{2}+1)\over \mu}v_{y}\right)u_{xy}
+\left(2v_{x}+{(\mu^{2}+1)\over \mu}u_{y}\right)v_{xy}
\nonumber \\
&+& \alpha_{y}-\beta_{y}.
\label{(5.28)}
\end{eqnarray}
The desired $\alpha$ and $\beta$ functions should satisfy $\alpha(\mu; x,y) < \beta(\mu;x,y)$, while their
gradient should guarantee that the partial derivatives of $H$ in (5.27) and (5.28) are never
simultaneously vanishing.

\section{Envelopes of our plane curve}

Our real manifold is a plane curve $\Gamma$ whose equation
\begin{equation}
H(\mu;x,y)=0
\label{(6.1)}
\end{equation}
involves an arbitrary parameter $\mu \in [0,1[$. As is well known, if each of the positions
of the curve $\Gamma$ is tangent to a fixed curve $E$, the curve $E$ is called the envelope of
the curves $\Gamma$, which are said to be enveloped by $E$.

If an envelope $E$ exists, let $(x,y)$ be the point of tangency of $E$ with that one
of the curves $\Gamma$ that corresponds to a certain value $\mu$ of the parameter. By
definition, the tangents to the curves $E$ and $\Gamma$ coincide for all values of $\mu$. 
If $\delta x$ and $\delta y$ are two quantities proportional to the direction cosines of the
tangent to $\Gamma$, and if ${dx \over d\mu}$ and ${dy \over d\mu}$ are the derivatives of
the unknown functions $x=\phi(\mu)$ and $y=\psi(\mu)$, the necessary condition for 
tangency is \cite{Goursat1904}
\begin{equation}
{{dx \over d\mu}\over \delta x}={{dy \over d\mu}\over \delta y}
\Longrightarrow {dy \over d\mu}\delta x-{dx \over d\mu}\delta y=0.
\label{(6.2)}
\end{equation}
On the other hand, since $\mu$ in Eq. (6.1) has a constant value for the particular
curve $\Gamma$ considered, we have
\begin{equation}
{\partial H \over \partial x}\delta x + {\partial H \over \partial y}\delta y=0,
\label{(6.3)}
\end{equation}
which determines the tangent to $\Gamma$. Moreover, the two unknown functions
$x=\phi(\mu),y=\psi(\mu)$ satisfy Eq. (6.1), where $\mu$ is now the independent variable.
Hence
\begin{equation}
{dH \over d \mu}={\partial H \over \partial x}{dx \over d \mu}
+{\partial H \over \partial y}{dy \over d\mu}+{\partial H \over \partial \mu}=0.
\label{(6.4)}
\end{equation}
The linear homogeneous system given by Eqs. (6.2) and (6.3) has nonvanishing solutions for
$\delta x$ and $\delta y$ if and only if the determinant of the $2 \times 2$ matrix of
coefficients vanishes, i.e.
\begin{equation}
{\partial H \over \partial x}{dx \over d \mu}+{\partial H \over \partial y}{dy \over d\mu}=0.
\label{(6.5)}
\end{equation}
By virtue of Eqs. (6.4) and (6.5) one finds that, if an envelope exists, its equation can be
found by eliminating the parameter $\mu$ between the equations $H=0$ and
\begin{equation}
{\partial H \over \partial \mu}=0.
\label{(6.6)}
\end{equation}

If $R(x,y)=0$ is the equation obtained by eliminating $\mu$ between (6.1) and (6.6), it can
be shown that it represents either the envelope of the curves $\Gamma$ or the locus of singular
points of these curves, at which
\begin{equation}
\left \{ H(\mu;x,y)=0, \; {\partial H \over \partial x}=0, \;
{\partial H \over \partial y}=0 \right \} \Longrightarrow 
{\partial H \over \partial \mu}=0.
\label{(6.7)}
\end{equation}
In other words, the plane curve $R(x,y)=0$ consists of two disjoint parts, one of which is the
envelope, while the other is the locus of singular points. 

When Eq. (6.1) is taken to have the
form (5.26), we find that the curve of equation $R(x,y)=0$ of the general theory outlined before
can be obtained, at least in principle, from
\begin{equation}
H=0 \Longrightarrow (u_{x})^{2}+(v_{x})^{2}+(u_{y})^{2}+(v_{y})^{2}
=\left(\mu+{1 \over \mu}\right)(u_{x}v_{y}-u_{y}v_{x})+\beta(\mu;x,y)-\alpha(\mu;x,y),
\label{(6.8)}
\end{equation}
\begin{equation}
{\partial H \over \partial \mu}=0 \Longrightarrow 
\left(1-{1 \over \mu^{2}}\right)(u_{x}v_{y}-u_{y}v_{x})
+{\partial \beta \over \partial \mu}-{\partial \alpha \over \partial \mu}=0.
\label{(6.9)}
\end{equation}
Interestingly, the construction of envelopes for the plane curve associated to pseudoholomorphic functions
of a complex variable leads to the elimination of the parameter $\mu \in [0,1[$ provided one is able
to solve the nonlinear partial differential equations (6.8) and (6.9).  

\subsection{The locus of points where $H,{\partial H \over \partial \mu},
{\partial^{2} H \over \partial \mu^{2}}$ all vanish}

Let us now consider a particular plane curve
\begin{equation}
H(\mu_{1};x,y)=0
\label{(6.10)}
\end{equation}
and a nearby curve, resulting from a small change of the parameter, i.e.
\begin{equation}
H(\mu_{1}+\lambda;x,y)=0.
\label{(6.11)}
\end{equation}
This equation can be replaced by
\begin{equation}
{H(\mu_{1}+\lambda;x,y)-H(\mu_{1};x,y) \over (\mu_{1}+\lambda)-\mu_{1}}=0,
\label{(6.12)}
\end{equation}
which converges, as $\lambda \rightarrow 0$, to the limit equation
\begin{equation}
\left({\partial H \over \partial \mu}\right)_{\mu=\mu_{1}}=0.
\label{(6.13)}
\end{equation}
Inspired by what is done in Ref. \cite{Bianchi1902} in the case of surfaces, we may therefore
consider the solution set of the three simultaneous equations 
\begin{equation}
H=0, \; {\partial H \over \partial \mu}=0, \;
H(\mu + \lambda;x,y)=0,
\label{(6.14)}
\end{equation}
the latter of which can be written in the form 
\begin{equation}
H+\lambda {\partial H \over \partial \mu}+{\lambda^{2}\over 2}{\partial^{2}H \over \partial \mu^{2}}
+\eta=0,
\label{(6.15)}
\end{equation}
where $\eta={\rm O}(\lambda^{3})$. By construction, this scheme is equivalent to studying the equations
\begin{equation}
H=0, \; {\partial H \over \partial \mu}=0, \;
{\partial^{2}H \over \partial \mu^{2}}+2{\eta \over \lambda^{2}}=0,
\label{(6.16)}
\end{equation}
where the latter equation converges, as $\lambda \rightarrow 0$, to the limit equation
\begin{equation}
{\partial^{2}H \over \partial \mu^{2}}=0.
\label{(6.17)}
\end{equation}

In our problem, our Eqs. (6.8) and (6.9) should be supplemented by
\begin{equation}
{\partial^{2}H \over \partial \mu^{2}}=0 \Longrightarrow
{2 \over \mu^{3}}(u_{x}v_{y}-u_{y}v_{x})
+{\partial^{2}\beta \over \partial \mu^{2}}-{\partial^{2}\alpha \over \partial \mu^{2}}=0.
\label{(6.18)}
\end{equation}
This scheme has interesting potentialities because, by virtue of (6.9), we can
write Eq. (6.18) in the form
\begin{equation}
{2 \over \mu(1-\mu^{2})}\left({\partial \beta \over \partial \mu}
-{\partial \alpha \over \partial \mu}\right)+{\partial^{2}\beta \over \partial \mu^{2}}
-{\partial^{2}\alpha \over \partial \mu^{2}}=0,
\label{(6.19)}
\end{equation}
which implies that, for a suitable unknown function $\sigma: x,y \rightarrow
\sigma(x,y)$, one can write
\begin{equation}
\gamma \equiv {\partial \beta \over \partial \mu}
-{\partial \alpha \over \partial \mu}
=\left({1 \over \mu^{2}}-1 \right)\sigma(x,y).
\label{(6.20)}
\end{equation}
Note that, in full agreement with what we say elsewhere in our paper, our formulas
are nontrivial only if $\mu \not=1$.    

\section{Concluding remarks}

In our paper, we have defined in Sect. $3$ the concept of pseudoholomorphic function
of $n$ complex variables. Moreover, in the case of one complex dimension, we have shown
in Sect. $5$ that the differential inequalities describing 
pseudoholomorphicity can be used to define a
one-real-dimensional manifold (by the vanishing of a function with nonzero gradient), 
which is here a $1$-parameter family of plane curves. In particular,
if the functions $u$ and $v$ satisfying these equations are taken
to be of the form (5.8) and (5.9), our construction becomes equivalent to obtaining a 
positive-definite quadratic form in the $(x,y)$ variables, for which (5.11)-(5.13), (5.20) 
and (5.21) should hold. On studying the envelopes associated to 
our plane curves, the parameter of 
the general theory can be eliminated by solving two nonlinear partial differential equations.
As far as we know, the consideration of such properties never appeared before in the literature,
although the basic concepts of our analysis were all well-known, when taken separately.

As far as we can see, no confusion should arise with the construction of pseudo-holomorphic submanifolds
performed in Ref. \cite{Taubes1996}, where the author starts instead from a compact, connected
$4$-manifold $X$ with a symplectic form $\omega$, and considers an almost complex structure $J$ for
the tangent bundle $TX$, $J$ being chosen to be compatible with the form $\omega$. Given a compact
submanifold $\Sigma$ of $X$, $\Sigma$ is said to be pseudo-holomorphic when $J$ maps the tangent
bundle $T \Sigma$ to itself as a subspace of $\left . TX \right |_{\Sigma}$.

It would be interesting to build one-complex-dimensional
manifolds whose transition functions are pseudoholomorphic 
according to the global theory of our Secs.  $2$ and $3$, but the differential inequalities
involved made it difficult for us to make progress along this line. 
In the present paper we have instead focused on one-real-dimensional manifolds
associated to the differential inequalities (2.23) of the global theory. 
Our next goal will be to get rid of the 
unknown parameter appearing in the definition of the real manifold by studying the geometry 
of plane curves, which will be reported in a future work \cite{future-ER}.

\acknowledgments
G. Esposito is grateful to the Dipartimento di Fisica ``Ettore Pancini'' of Federico II University for
hospitality and support. The research of R. Roychowdhury was supported by FAPESP through 
Instituto de Fisica, Universidade de Sao Paulo with grant number 2013/17765-0. This work was 
initiated during RR's visit to Federico II University in Naples. He thanks INFN Napoli 
for the hospitality and support during that period.

\begin{appendix}

\section{Counterexamples}

Suppose that, for simplicity, we try to reduce the differential inequalities of Secs. $2$ and $3$
to equalities. For the two-variable case, for $k=1$ 
if we saturate the first of the  inequalities (\ref{(3.23)}) 
we obtain two independent equalities 
\begin{equation}
\mu_1 J_1 = \varphi_1^{2}, \;\;\;  \Phi_1^{2} = {J_1 \over \mu_1}
\label{(A1)}
\end{equation}

Now  from (\ref{(A1)}) one can see that the ratio of the squares of upper and lower limits of the increment ratios
\begin{equation}
\frac{\Phi_1^{2}}{\varphi_1^{2}} = \frac{1}{\mu_1^{2}},
\label{(A2)}
\end{equation}
which in turn yields for $\mu_1 < 1$
\begin{equation}
\Phi_1^{2} > \varphi_1^{2},
\label{(A3)}
\end{equation}
and this minorization is consistent with our definitions.

We note that we can form a nonlinear partial differential equation  involving  
$u_{x_1},u_{y_1},v_{x_1},v_{y_1}$ by using equations 
(\ref{(3.21)}) and (\ref{(A1)}) and writing
\begin{equation}
\mu_1 J_1 = \frac{1}{2}(E_1+G_1- \omega_1)
\label{(A4)}
\end{equation}
where
\begin{equation}
E_1 \equiv (u_{x_1})^{2}+(v_{x_1})^{2}, \;
G_1 \equiv (u_{y_1})^{2}+(v_{y_1})^{2}, \;
F_1 \equiv u_{x_1}u_{y_1}+v_{x_1}v_{y_1},
\label{(A5)}
\end{equation}
and
\begin{equation}
\omega_1 \equiv \sqrt{(E_1-G_1)^{2}+4F_1^{2}}.
\label{(A6)}
\end{equation}
Since $J_1 \neq 0$, after a little algebra one can find the following differential equation:
\begin{equation}
(\mu_1 + 1)(u_{x_1}v_{y_1}-u_{y_1}v_{x_1}) - [(u_{x_1})^2 + (u_{y_1})^2 
+ (v_{x_1})^2 +( v_{y_1})^2]\mu_1 = 0,
\label{(A7)}
\end{equation}
where we want to solve for $u(x_1,y_1)$ and $v(x_1,y_1)$.
 
From here onwards we will be writing $x_1 = x$ and $y_1 = y$ and  the parameter 
$\alpha = \alpha_{\mu_1} = 1+\frac{1}{\mu_1}$
so that we settle for this partial differential equation
\begin{equation}
\alpha (u_{x}v_{y}-u_{y}v_{x}) = (u_{x})^2 + (u_{y})^2 + (v_{x})^2 +( v_{y})^{2}. 
\label{(A8)}
\end{equation}

Now for $\mu_1 = 1$ i.e. $\alpha=2$ we get the class of holomorphic functions with Cauchy-Riemann 
conditions being satisfied, since in this case
one can write (\ref{(A7)}) as 
$$
(u_{x}-v_{y})^2 +( u_{y}+v_{x})^2 =0.
$$

In order to get the most general solution for the class of pseudo analytic functions we need to solve for
\begin{equation}
(\alpha - 2) (u_{x}v_{y}-u_{y}v_{x}) = (u_{x}-v_{y})^2 +( u_{y}+v_{x})^{2}. 
\label{(A9)}
\end{equation}
As $\mu_1 = \frac{1}{2}$ and hence $\alpha = 3$ is an admissible choice we might try to find a solution 
for the following nonlinear PDE in 2 variables with constant coefficients
\begin{equation}
(u_{x}v_{y}-u_{y}v_{x}) = (u_{x}-v_{y})^2 +( u_{y}+v_{x})^{2}. 
\label{(A10)}
\end{equation}
But even before solving this, if we try to saturate the other inequalitiy in (\ref{(3.23)}) 
and work out the case $\mu_1 = \frac{1}{2}$ we soon run into inconsistencies as we 
end up getting an equation which is 
\begin{equation}
\frac{1}{2} (u_{x}v_{y}-u_{y}v_{x}) = (u_{x}-v_{y})^2 +( u_{y}+v_{x})^{2}. 
\label{(A11)}
\end{equation}
This clearly hints at possible inconsistencies that we head towards with our attempt of 
saturating the set of inequalities. The conclusion 
is that we need strict inequalities in order to deal with pseudoholomorphic functions.

As a second example, let $w$ be the holomorphic function 
\begin{equation}
w(z)= z^2 = u(x,y)+{\rm i}v(x,y)
\label{(A12)}
\end{equation}
so that, by definition, the functions $u(x,y),v(x,y)$ read as 
\begin{equation}
 u(x,y) = (x^2 - y^2), \;\;\;\; v(x,y) = 2xy.
\label{(A13)}
\end{equation}
Now let us try to deform $w(z)$ by adding to it a small non-holomorphic part $\epsilon \bar z$ 
such that the modified functions $u(x,y),v(x,y)$ become 
\begin{equation}
u(x,y) = (x^2 - y^2 + \epsilon x), \;\;\;\; v(x,y) = (2xy - \epsilon y).
\label{(A14)}
\end{equation}
With this it is easy to see that the Cauchy-Riemann condition breaks down as 
$u_{x} \neq v_{y}$  although  $u_{y}=-v_{x}$. This is clearly a non-analytic case,
but can we go ahead with this and recover, 
with a suitable choice of $\mu <1$, the pseudoholomorphic case?

Let us present a few computational details.
The Jacobian defined in (\ref{(2.2)}) for this case becomes 
\begin{equation}
J (x,y) = 4(x^2 + y^2) - \epsilon^{2}.
\label{(A15)}
\end{equation}
If we make use of the  master inequality (\ref{(2.5)}) we find that, 
for a positive real number  $\mu \in ]0,1]$, it is always true that
\begin{equation}
\varphi^2(x,y) \geq \mu^2 \Phi^2(x,y).
\label{(A16)}
\end{equation}
Using the identitites (\ref{(2.21)}) we get a ratio satisfying the following bound $\forall \mu$:
\begin{equation}
\frac {E+G- \omega}{E+G+ \omega} \geq \mu^{2}, 
\label{(A17)}
\end{equation}
with $E$, $G$ and $\omega$ respectively given by
\begin{eqnarray}
E  =  (u_{x})^{2}+(v_{x})^{2} = (2x + \epsilon)^2 + 4y^{2}, \nonumber\\
G = (u_{y})^{2}+(v_{y})^{2} = 4y^2 +  (2x - \epsilon)^{2}, \nonumber\\
\omega = \sqrt{(E+G)^{2}-4J^{2}} = 8 \sqrt{(x^2 + y^2)}\epsilon.
\label{(A18)}
\end{eqnarray}
Upon using (\ref{(A16)}) and (\ref{(A17)}) , $\forall \epsilon$ small we get the following bound on $\epsilon$:
$$
2(x^2 + y^2) > \epsilon^2.
$$
However, since $\epsilon$ is very small but finite, as $(x,y)$ tends to $(0,0)$ this condition is violated 
on a set of finite measure because the Lebesgue measure  of the set 
$$A_{\epsilon}=
\left\{x,y : 0 \leq 2(x^2 + y^2) < \epsilon^{2} \right\}
$$
is nonvanishing.
But then the condition for pseudoholomorphicity would not be satisfied almost everywhere. This means that 
our example fails to provide a pseudoholomorphic function.

\end{appendix}


\begin{thebibliography}{}
\bibitem{Penrose1975}
R. Penrose, Twistor theory, its aims and achievements, in 
{\it Quantum Gravity, an Oxford Syposium}, eds. C. J. Isham, R. Penrose and D. W. Sciama
(Clarendon Press, Oxford, 1975).
\bibitem{gesposit-complex}
G. Esposito,  {\it Complex General Relativity} (Kluwer Academic Publishers, Dordrecht, 1995).
\bibitem{Penrose1976}
R. Penrose, Non-linear gravitons and curved twistor theory, {\it Gen. Rel. Grav.} 
{\bf 7} (1976) 31.
\bibitem{Penrose1986}
R. Penrose, Twistors in general relativity, in {\it General Relativity and Gravitation}, ed.
M. A. H. MacCallum (Cambridge University Press, Cambridge, 1986).
\bibitem{Chern1979}
S. S. Chern, {\it Complex Manifolds without Potential Theory} (Springer, Berlin, 1979).
\bibitem{Goursat1884}
E. Goursat, Proof of the Cauchy theorem, {\it Acta Math.} {\bf 4} (1884) 197.
\bibitem{Goursat1900}
E. Goursat, On the general definition of analytic functions, following Cauchy,
{\it Trans. Amer. Math. Soc.} {\bf 1} (1900) 14.
\bibitem{Bers1956}
L. Bers, An outline of the theory of pseudoanalytic functions, 
{\it Bull. Amer. Math. Soc.} {\bf 62} (1956) 291.
\bibitem{Caccioppoli1952}
R. Caccioppoli, Foundations for a general theory of pseudoanalytic functions of a complex
variable, {\it Rend. Acc. Naz. Lincei}, Ser. VIII, {\bf 13} (1952) 197; 
ibid. {\bf 13} (1952) 321.
\bibitem{Caccioppoli1953}
R. Caccioppoli, Pseudo-analytic functions and pseudo-conformal representations of
Riemann surfaces, {\it Ricerche Mat.} {\bf 2} (1953) 104.
\bibitem{Courant1962}
R. Courant and D. Hilbert, {\it Methods of Mathematical Physics, Vol. II: Partial 
Differential Equations} (Wiley, New York, 1962).
\bibitem{Hedrick1933}
E. R. Hedrick, Non-analytic functions of a complex variable, 
{\it Bull. Amer. Math. Soc.} {\bf 39} (1933) 75.
\bibitem{Donaldson1989}
S. K. Donaldson and D. P. Sullivan, Quasiconformal 4-manifolds,
{\it Acta Math.} {\bf 163} (1989) 181.
\bibitem{Bers1957}
L. Bers, On a theorem of Mori and on the definition of quasi-conformality,
{\it Trans. Amer. Math. Soc.} {\bf 84} (1957) 78.
\bibitem{Morrey1938}
C. B. Morrey, On the solutions of quasi-linear elliptic partial differential equations,
{\it Trans. Amer. Math. Soc.} {\bf 43} (1938) 126.
\bibitem{Mori1957}
A. Mori, On quasi-conformality and pseudo-analyticity, 
{\it Trans. Amer. Math. Soc.} {\bf 84} (1957) 56.
\bibitem{Hitotumatu1959}
S. Hitotumatu, On quasi-conformal functions of several complex variables,
{\it J. Math. Mech.} {\bf 8} (1959) 77.
\bibitem{Toki1954}
Y. Toki and K. Shibata, On the pseudo-analytic functions,
{\it Osaka Math. J.} {\bf 6} (1954) 145.
\bibitem{Storvick1957}
D. A. Storvick, On pseudo-analytic functions, {\it Nagoya Math. J.} {\bf 12} (1957) 131.
\bibitem{Koohara1971}
A. Koohara, Similarity principle of the generalized Cauchy-Riemann equations for several
complex variables, {\it J. Math. Soc. Japan} {\bf 23} (1971) 213.
\bibitem{Fryant1981}
A. J. Fryant, Ultraspherical expansions and pseudoanalytic functions, 
{\it Pac. J. Math.} {\bf 94} (1981) 83.
\bibitem{Tutschke2007}
W. Tutschke, Generalized analytic functions in higher dimensions, 
{\it Georgian Math. J.} {\bf 14} (2007) 581.
\bibitem{Rheinboldt1996}
W. C. Rheinboldt, MANPAK: A set of algorithms for computations on implicitly defined manifolds,
{\it Comp. Math. Appl.} {\bf 32} (1996) 15.
\bibitem{Schwartz1967}
L. Schwartz, {\it Cours d'Analyse} (Hermann, Paris, 1967).
\bibitem{Goursat1904}
E. Goursat, {\it A Course in Mathematical Analysis, Vol. 1} (Ginn and Company, Boston, 1904).
\bibitem{Bianchi1902}
L. Bianchi, {\it Lectures on Differential Geometry} (Enrico Spoerri, Pisa, 1902).
\bibitem{Taubes1996}
C. H. Taubes, Counting pseudo-holomorphic submanifolds in dimension $4$, 
{\it J. Diff. Geom.} {\bf 44} (1996) 818.
\bibitem{future-ER}
G. Esposito and R. Roychowdhury, in progress.
\end{thebibliography}
\end{document}